\documentclass{amsart}
\usepackage{amsfonts, amsmath, amsthm, amscd, amsfonts, amssymb, graphicx, color, xypic, IMjournal}
\usepackage{amsfonts}
\usepackage{mathrsfs}
\begin{document}
\newtheorem{The}{Theorem}[section]
\newtheorem{Def}[The]{Definition}
\newtheorem{Cor}[The]{Corollary}
\newtheorem{Lem}[The]{Lemma}
\newtheorem{Prop}[The]{Proposition}
\newtheorem{Rem}[The]{Remark}
\newtheorem{Exa}[The]{Example}
\newtheorem{Claim}[The]{Claim}
\def\pn{\par\noindent}
\def\cen{\centerline}

\numberwithin{equation}{section}

\title{Gorenstein projective complexes with respect to cotorsion pairs}

\author{\|Renyu |Zhao|, \|Pengju |Ma|, Lanzhou}

%\rec {January 27, 2006}

\dedicatory{Cordially dedicated to ...}

\abstract
%%%%%%
%%%7 to 12 lines optimum; references in the abstract should be given in full form, e.g. V. Nov\'ak, M. Novotn\'y: Linear extensions of orderings. Czech. Math. J. 50 (2000), 853--864.
%%%%%
Let $(\mathcal{A,B})$ be a complete and hereditary cotorsion pair in the category of
left $R$-modules. In this paper, the so-called Gorenstein projective complexes respect to the cotorsion pair $(\mathcal{A}, \mathcal{B})$ are introduced. We show that these complexes are just the complexes of Gorenstein projective modules respect to the cotorsion pair $(\mathcal{A}, \mathcal{B})$. As applications, we prove that both Gorenstein projective modules with respect to cotorsion pairs and Gorenstein projective complexes with respect to cotorsion pairs possess of stability.
\endabstract



\thanks
   The research was partially supported by the National Natural Science Foundation of China (Nos. 11361051, 11361052, 11401476)
\endthanks

\section{Introduction}\label{sec1}

%%%%%
Let $(\mathcal{A,B})$ be a complete and hereditary cotorsion pair in the category of
left $R$-modules. Then there is two induced cotorsion pairs $(\mathcal{\widetilde{A},\text{dg}\widetilde{B}})$ and $(\mathcal{\text{dg}\widetilde{A},\widetilde{B}})$ (\cite{g04}), and both of them are complete and hereditary (\cite{yd,yl11}). Recently, among others, the Gorenstein category $\mathcal{G}(\mathcal{A})$ with respect to the cotorsion pair $(\mathcal{A,B})$ was introduced and studied by Yang and Chen \cite{yc} (see Definition \ref {def}). In this paper, we generalize this notion to the category of complexes of left $R$-modules, namely, we introduce the Gorenstein projective complexes respect to the cotorsion pair $(\mathcal{A}, \mathcal{B})$ (see Definition \ref {D1}). The class of these complexes will be denoted by $\mathcal{G}(\widetilde{\mathcal{A}})$. It contains Gorenstein projective complexes \cite{eg}, \textbf{F}-Gorenstein flat complexes \cite{hx}, and Gorenstein flat complexes \cite{gr} over right coherent rings as its special cases. By using the techniques of Bravo and Gillespie \cite{bg}, we prove the following result (see Theorem \ref{T1}).

\begin{The}\label{T1} A complex $C$ of left $R$-modules belongs to $\mathcal{G}(\mathcal{\widetilde{A}})$ if and only if each $C_i\in\mathcal{G}(\mathcal{A})$.\end{The}

This result unifies and generalizes of \cite[Theorems 4.7]{hx} and \cite[Theorems 2.2, 3.1]{yl}.
By this connection, we prove the following results (see Theorem \ref{T2} and Theorem \ref{T3}, respectively).

\begin{The} A left $R$-module $M$ belongs to $\mathcal{G}(\mathcal{A})$ if and only if there exists a ${\rm Hom}_{R}(-,\mathcal{A}\cap \mathcal{B})$-exact exact sequence $\xymatrix@-1pc{\cdots\ar[r]^{} & G_{1} \ar[r]& G_{0}\ar[r] &G_{-1}\ar[r] &\cdots}$
in $\mathcal{G}(\mathcal{A})$ such that $M\cong{\rm Im}(\xymatrix@-1pc{G^{0}\ar[r] &G^{-1}})$.
\end{The}

\begin{The}\label{T1} A complex $C$ of left $R$-modules belongs to $\mathcal{G}(\mathcal{\widetilde{A}})$ if and only if there exists a ${\rm Hom}_{C(R)}(-,\mathcal{\widetilde{A}}\cap \text{dg}\mathcal{\widetilde{B}})$-exact exact sequence $\xymatrix@-1pc{\cdots\ar[r]^{} & G^{1} \ar[r]& G^{0}\ar[r] &G^{-1}\ar[r] &\cdots}$
in $\mathcal{G}(\mathcal{\widetilde{A}})$ such that $C\cong{\rm Im}(\xymatrix@-1pc{G^{0}\ar[r] &G^{-1}})$.\end{The}

These two results imply that the category $\mathcal{G}(\mathcal{A})$ and  $\mathcal{G}(\mathcal{\widetilde{A}})$ have stability, respectively.

The contents of this paper are summarized as follows. In Section 2, we
review some basic notation and notions for use throughout this paper. Section 3 is devoted to introduce the notion of Gorenstein projective complexes respect to cotorsion pairs and give the proof Theorem 1.1. By using of Theorem 1.1, in Section 4, we give the proof of Theorem 1.2 and Theorem 1.3.

\section{Preliminaries}\label{sec1}

Throughout this article, $R$ denotes an associative ring with
identity, modules are assumed to be unitary, and the default action of the ring is on the
left. Right modules over $R$ are hence treated as (left) modules over the opposite ring
$R^{\circ}$. We use $R$-Mod to denote the category of $R$-modules, $C(R)$ to denote the category of complexes of $R$-modules and $\mathcal{P}$ (resp., $\mathcal{F}, \mathcal{C}$) to denote the class of projective (resp., flat, cotorsion) $R$-modules.

A complex $$\xymatrix{\cdots \ar[r]^{} &C_{n+1} \ar[r]^{\delta_{n+1}} &  C_{n}  \ar[r]^{\delta_{n}}& C_{n-1} \ar[r]^{}&\cdots}$$
will be denoted by $(C,\delta)$ or simply $C$. The $n$th cycle (resp. boundary, homology) of
$C$ is denoted by ${\rm Z}_n(C)$ (resp., ${\rm B}_n(C)$, ${\rm H}_n(C)$). We will use superscripts to distinguish complexes. So if $\{C^i\}_{i\in I}$ is a family of complexes, $C^i$ will be complex $$\xymatrix{\cdots \ar[r]^{} &C_{n+1}^i \ar[r]^{\delta_{n+1}^i} &  C_{n}^i  \ar[r]^{\delta_{n}^i}& C_{n-1}^i \ar[r]^{}&\cdots.}$$ Given an $R$-module $M$, we will denote by $\overline{M}$ the complex
$$\xymatrix{\cdots \ar[r]^{} &0 \ar[r]^{} & M \ar[r]^{id}& M\ar[r]^{}&0 \ar[r]^{}&\cdots}$$
with $M$ in 1 and 0th degrees. Given a $C\in C(R)$ and an integer $m$,
$C[m]$ denotes the complex such that $C[m]_n=C_{n-m}$ and whose boundary operators
are $(-1)^m\delta_{n-m}$. Given $C,D\in C(R)$, we use $\text{Hom}_{C(R)}(C,D)$ to present the group of all morphisms
from $C$ to $D$, and Ext$^{i}_{C(R)}(C,D)$ denotes the groups one gets from the right derived functor of Hom for $i\geq 0$.

Let $C\in C(R^{\circ})$ and $D\in C(R)$.
The tensor product $C\otimes_R D$ is the $\mathbb{Z}$-complex whose underlying graded module is
given by $(C\otimes_R D)_n=\bigoplus_{i+j=n}C_i\otimes_R D_j$, and whose differential is defined by specifying its action on an elementary tensor
of homogeneous elements as $\delta^{C\otimes_R D}(x\otimes y)=\delta^C(x)\otimes y + (-1)^{|x|}
x\otimes \delta^D(y)$, where $|x|$ is the degree of $x$ in $C$. Let $C\overline{\otimes}_R D =\frac{C\otimes_R D}{{\rm B}(C\otimes_R D)}$
, that is, $C\overline{\otimes}_R D$ is the complex of abelian groups
with $n$th entry $(C\overline{\otimes}_R D)_n=\frac{(C\otimes_R D)_n}{{\rm B}_n(C\otimes_R D)}$ and boundary map $\delta^{C\overline{\otimes}_R D}(\overline{x\otimes y})=\overline{\delta^C(x)\otimes y}$, where $\overline{x\otimes y}$ is used to denote the coset in $\frac{(C\otimes_R D)_n}{{\rm B}_n(C\otimes_R D)}$.
This gives us new right exact bifunctor $-\overline{\otimes}_R -$ which have left
derived functor $ \overline{{\rm Tor}}_i(-,-)$.

For  $C,D\in C(R)$, ${\rm Hom}_R(C,D)$ is the complex of abelian groups with the degree-$n$ term $\text{Hom}_R(C,D)_n=\prod_{i\in \mathbb{Z}}\text{Hom}_R(C_i,D_{n+i}),$ and its boundary operators are $\delta_n^{{\rm Hom}_R(C,D)}\left((f_i)_{i\in \mathbb{Z}}\right)=\left(\delta^{D}_{n+i}f_i-{(-1)^n} f_{i-1}\delta^{C}_{i}\right)_{i\in \mathbb{Z}}$ for any $(f_i)_{i\in \mathbb{Z}}\in
\text{Hom}_R(C,D)_n$.
Let ${\rm \underline{Hom}}_R(C,D)= {\rm Z}({\rm Hom}_R(C,D))$, that
is, ${\rm \underline{Hom}}_R(C,D)$ is the complex of abelian groups with $n$th entry ${\rm \underline{Hom}}_R(C,D)_n
={\rm Z}_n({\rm Hom}_R(C,D))={\rm Hom}_{C(R)}(C,D[-n])$ and boundary map $\delta_n^{{\rm \underline{Hom}}_R(C,D)}\left((f_i)_{i\in \mathbb{Z}}\right)=\left((-1)^n\delta^{D}_{n+i}f_i\right)_{i\in \mathbb{Z}}$ for any $(f_i)_{i\in \mathbb{Z}}\in
\text{Hom}_R(C,D)_n$. Then we get new functors ${\rm \underline{Hom}}_R(C,-)$ and ${\rm \underline{Hom}}_R(-,D)$ which
are both left exact. The book \cite{gr} is a standard reference for complexes.

Let $\mathcal{D}$ be an abelian category. A pair $(\mathcal{A}, \mathcal{B})$ of
classes of objects of $\mathcal{D}$ is called a cotorsion pair if
$\mathcal{A}^{\bot}=\mathcal{B}$ and $\mathcal{A }={^{\bot}\mathcal{B}}$, where $\mathcal{A}^{\bot}=\{D\in \mathcal{D}\mid{\rm Ext}_{\mathcal{D}}^1(A,D)=0 \text{ for all } A\in \mathcal{A}\}$ and ${^{\bot}\mathcal{B}}=\{D\in \mathcal{D}\mid{\rm Ext}_{\mathcal{D}}^1(D,B)=0 \text{ for all } B\in \mathcal{B}\}$.
%An abelian category $\mathcal{D}$ is said to have enough $\mathcal{A}$-objects (resp., enough $\mathcal{B}$-objects) if for every object $D\in \mathcal{D}$ there is an exact sequence $\xymatrix@-1pc{A \ar[r] &D \ar[r]& 0}$ (resp., $\xymatrix@-1pc{0 \ar[r] &D \ar[r]& B}$), where $A\in \mathcal{A}$ (resp., $B\in \mathcal{B}$).
A special $\mathcal{A}$-precover (resp., special $\mathcal{B}$-preenvelope) of an object $D\in \mathcal{D}$ is a
short exact sequence $\xymatrix@-1pc{0 \ar[r]& B \ar[r]& A \ar[r] &D \ar[r]& 0}$ (resp., $\xymatrix@-1pc{0 \ar[r]& D \ar[r]& B' \ar[r]& A' \ar[r]& 0}$), where
$A\in \mathcal{A}$ and $B\in \mathcal{B}$ (resp., $A'\in \mathcal{A}$ and $B'\in \mathcal{B}$). A cotorsion pair $(\mathcal{A}, \mathcal{B})$ is said to be complete if every
object $D\in \mathcal{D}$ has a special $\mathcal{A}$-precover and a special $\mathcal{B}$-preenvelope.
A cotorsion pair $(\mathcal{A}, \mathcal{B})$ in $\mathcal{D}$ is said to be hereditary if ${\rm Ext}^i_{\mathcal{D}}(A,B)=0$ for all $A\in \mathcal{A},B\in \mathcal{B}$ and all $i\geqslant 1$.
%The class $\mathcal{A}$ is called Projective resolving if $\mathcal{A}$ is closed under kernels of epimorphisms and the class $\mathcal{B}$ is called injective coresolving if $\mathcal{B}$ is closed under cokernels of monomorphisms. Note that for a hereditary cotorsion pair $(\mathcal{A}, \mathcal{B})$, the classes $\mathcal{A}$ and $\mathcal{B}$ are resolving and coresolving, respectively. Conversely, if $\mathcal{D}$ has enough $\mathcal{A}$-objects and enough $\mathcal{B}$-objects, then $(\mathcal{A}, \mathcal{B})$ is a hereditary cotorsion pair if and only if $\mathcal{A}$ is a resolving class if and only if $\mathcal{B}$ is a coresolving class.
If we choose $\mathcal{D}=R\text{-Mod}$ for some ring $R$, the most obvious example of a complete hereditary cotorsion pair is $(\mathcal{P}, R\text{-Mod}).$ Perhaps one of the most useful complete hereditary cotorsion pair is the flat cotorsion pair $(\mathcal{F},\mathcal{ C})$. For a good reference on cotorsion pairs see \cite{ejbook}.

\begin{Def} {\rm(\cite[Definition 3.3]{g04})}
   Let $(\mathcal{A}, \mathcal{B})$ be a cotorsion pair in $R$-Mod
and $X$ an $R$-complex.\\
\indent(1) $X$ is called an $\mathcal{A}$ complex if it is exact and ${\rm Z}_n(X)\in \mathcal{A}$ for all $n\in \mathbb{Z}$.\\
\indent(2) $X$ is called a $\mathcal{B}$ complex if it is exact and ${\rm Z}_n(X)\in \mathcal{B}$ for all $n\in \mathbb{Z}$.\\
\indent(3) $X$ is called a dg-$\mathcal{A}$ complex if each $X_n\in \mathcal{A}$ and ${\rm Hom}_R(X,B)$ is
exact whenever $B$ is a $\mathcal{B}$ complex.\\
\indent(4) $X$ is called a dg-$\mathcal{B}$ complex if each $X_n\in \mathcal{B}$ and ${\rm Hom}_R(A,X)$ is
exact whenever $A$ is an $\mathcal{A}$ complex.\end{Def}
\indent We denote the class of $\mathcal{A}$ complexes by $\widetilde{\mathcal{A}}$ and the class of dg-$\mathcal{A}$ complexes by dg$\widetilde{\mathcal{A}}$.
Similarly, the $\mathcal{B}$ complexes are denoted by $\widetilde{\mathcal{B}}$ and the class of dg-$\mathcal{B}$ complexes are
denoted by dg$\widetilde{\mathcal{B}}$.  We sometimes name $\mathcal{A}$ (resp., dg-$\mathcal{A}$) complexes by the name of the class $\mathcal{A}$. For example, the projective (resp. dg-projective) complexes are actually the $\mathcal{P}$ (resp.
dg-$\mathcal{P}$) complexes. It follows from \cite[Proposition 3.6]{g04} that $(\widetilde{\mathcal{A}}, \text{dg}\widetilde{\mathcal{B}})$ and $(\text{dg}\widetilde{\mathcal{A}}, \widetilde{\mathcal{B}})$ are cotorsion pairs in $C(R)$. Moreover, by \cite[Corollary 3.13]{g04}, \cite[Theorem 2.4 and Corollary 2.7]{yd} or \cite[Theorem 3.5]{yl11}, we have the following facts.

\begin{Lem} Let $(\mathcal{A}, \mathcal{B})$ be a complete hereditary cotorsion pair in $R\text{-Mod}$. Then the induced cotorsion pairs $(\widetilde{\mathcal{A}}, \text{dg}\widetilde{\mathcal{B}})$ and $(\text{dg}\widetilde{\mathcal{A}}, \widetilde{\mathcal{B}})$ in $C(R)$ are both complete and hereditary. Furthermore, $\text{dg}\widetilde{\mathcal{A}}\cap\mathcal{E}=\widetilde{\mathcal{A}}$ and $\text{dg}\widetilde{\mathcal{B}}\cap \mathcal{E}=\widetilde{\mathcal{B}},$ where $\mathcal{E}$ is the class of exact complexes.
\end{Lem}

Let $\mathcal{D}$ be an abelian category and $\mathcal{H}$ a full subcategory of $\mathcal{D}$. Recall that a sequence $\textbf{S}$ in $\mathcal{D}$ is ${\rm Hom}_{\mathcal{D}}(-,\mathcal{H})$-exact (resp., ${\rm Hom}_{\mathcal{D}}(\mathcal{H},-)$-exact) if the sequence ${\rm Hom}_{\mathcal{D}}(\textbf{S},H)$ (resp., ${\rm Hom}_{\mathcal{D}}(H,\textbf{S})$) is exact for any $H\in \mathcal{H}$.

\begin{Def}\label{def} {\rm(\cite[Definition 3.1]{yc})} Let $(\mathcal{A}, \mathcal{B})$ be a complete
and hereditary cotorsion pair in $R$-Mod. An $R$-module $M$ is called Gorenstein projective respect to the cotorsion pair $(\mathcal{A}, \mathcal{B})$ if
there exists a ${\rm Hom}_{R}(-,\mathcal{A}\cap\mathcal{B})$-exact exact sequence $\xymatrix@-0.8pc{\cdots\ar[r]^{} & A_{1} \ar[r]& A_{0}\ar[r] &A_{-1}\ar[r] &\cdots}$ with each $A_i\in\mathcal{A}$, such that $M\cong{\rm Im}(\xymatrix@-1pc{A_{0}\ar[r] &A_{-1}})$. We let $\mathcal{G}(\mathcal{A})$ be the class of Gorenstein projective $R$-modules respect to the cotorsion pair $(\mathcal{A}, \mathcal{B})$.\end{Def}

\begin {Rem}\label{rem1} (1) By completeness of the cotorsion pair $(\mathcal{A}, \mathcal{B})$, an $R$-module $M$ is in $\mathcal{G}(\mathcal{A})$ if and only if ${\rm Ext}^{i\geqslant 1}_R(M,N)=0$ for any $N\in\mathcal{A}\cap\mathcal{B}$ and there exists a ${\rm Hom}_R(-,\mathcal{A}\cap\mathcal{B})$-exact exact sequence $\xymatrix@-1pc{0\ar[r]^{} & M \ar[r]& A_{0}\ar[r] &A_{-1}\ar[r] &\cdots}$ with each $A_i\in \mathcal{A}$.
\par(2) This definition unifies the following notions: Gorenstein projective modules \cite{ej,h} (in the case $(\mathcal{A}, \mathcal{B})=(\mathcal{P},R\text{-Mod})$); {\bf F}-Gorenstein flat modules \cite{hx} (when $(\mathcal{A}, \mathcal{B})=(\mathcal{F},\mathcal{C})$); and Gorenstein flat modules \cite{gr} (when $(\mathcal{A}, \mathcal{B})=(\mathcal{F},\mathcal{C})$ and $R$ is a right coherent ring), see \cite[Lemma 3.2]{hx}.

\end{Rem}

In what follows, we always assume that $(\mathcal{A}, \mathcal{B})$ is a complete
and hereditary cotorsion pair in $R$-Mod.

\section{Gorenstein projective complexes with respect to cotorsion pairs}\label{sec1}

\begin{Def}\label{D1} An $R$-complex $C$ is called Gorenstein projective respect to the cotorsion pair $(\mathcal{A}, \mathcal{B})$ if
there exists a ${\rm Hom}_{C(R)}(-,\mathcal{\widetilde{A}\cap \text{dg}\widetilde{B}})$-exact exact sequence $\xymatrix@-1pc{\cdots\ar[r]^{} & A^{1} \ar[r]& A^{0}\ar[r] &A^{-1}\ar[r] &\cdots}$ with each $A^i\in \mathcal{\widetilde{A}}$ such that $C\cong{\rm Im}(\xymatrix@-1pc{A^{0}\ar[r] &A^{-1}})$.\end{Def}
We denote the class of Gorenstein projective $R$-complexes respect to the cotorsion pair $(\mathcal{A}, \mathcal{B})$ by $\mathcal{G}(\widetilde{\mathcal{A}})$.

\begin {Rem}\label{example} (1) It is clear that $\widetilde{\mathcal{A}}\subseteq\mathcal{G}(\widetilde{\mathcal{A}})$. If $\xymatrix@-1pc{\textbf{A}=\cdots\ar[r]^{} & A^{1} \ar[r]& A^{0}\ar[r] &A^{-1}\ar[r] &\cdots}$ is a ${\rm Hom}_{C(R)}(-,\mathcal{\widetilde{A}\cap \text{dg}\widetilde{B}})$-exact exact sequence of complexes in $\mathcal{\widetilde{A}}$, then
by symmetry, all the images, the kernels and the cokernels of $\textbf{A}$ are in $\mathcal{G}(\widetilde{\mathcal{A}})$.

(2) If $(\mathcal{A}, \mathcal{B})=(\mathcal{P}, R\text{-\rm Mod})$, then Gorenstein projective complexes respect to the cotorsion pair $(\mathcal{A}, \mathcal{B})$ are exactly the Gorenstein projective complexes in \cite{eg}.

(3) If $(\mathcal{A}, \mathcal{B})=(\mathcal{F}, \mathcal{C})$, the flat cotorsion pair, then Gorenstein projective complexes respect to the cotorsion pair $(\mathcal{A}, \mathcal{B})$ are just \textbf{F}-Gorenstein flat complexes in \cite{hx}.
\end{Rem}

Recall from \cite{gr} that a short exact sequence $\xymatrix@-1pc{0\ar[r]^{} & S \ar[r]& C\ar[r] &C/S\ar[r] &0}$ in $C(R)$ is pure if the sequence $\xymatrix@-1pc{0\ar[r]^{} & D\overline{\otimes}_R S \ar[r]& D\overline{\otimes}_R C \ar[r] &D\overline{\otimes}_R C/S \ar[r] &0}$ is exact for any $D\in C(R^{\circ})$. According to \cite{gr}, an $R$-complex $C$ is called Gorenstein flat if there exists an
exact sequence of flat complexes $\xymatrix@-1pc{\cdots\ar[r]^{} & F^{1} \ar[r]& F^{0}\ar[r] &F^{-1}\ar[r] &\cdots}$ with
$C\cong{\rm Im}(\xymatrix@-1pc{F^{0}\ar[r] &F^{-1}})$ and which remains exact after applying $I\overline{\otimes}_R-$ for any injective
$R^{\circ}$-complex $I$. The next result shows that Gorenstein projective complexes respect to the cotorsion pair $(\mathcal{F}, \mathcal{C})$ over right coherent rings are just Gorenstein flat complexes.

\begin{Prop}\label{prop}  If $R$ is a right coherent ring, then $C$ is an \textbf{F}-Gorenstein flat complexes if and only if $C$ is Gorenstein flat.\end{Prop}

\proof $\Longrightarrow$) Assume that $C$ is an \textbf{F}-Gorenstein flat complex. Then there exists a ${\rm Hom}_{C(R)}(-,\mathcal{\widetilde{F}\cap \text{dg}\widetilde{C}})$-exact exact sequence $\xymatrix@-1pc{\cdots\ar[r]^{} & F^{1} \ar[r]& F^{0}\ar[r] &F^{-1}\ar[r] &\cdots}$ of flat complexes such that $C\cong{\rm Im}(\xymatrix@-1pc{F^{0}\ar[r] &F^{-1}})$. Let $I$ be an injective complex of right $R$-modules. Then $I^{+}[-n] \in \mathcal{\widetilde{F}\cap \text{dg}\widetilde{C}}$ for any $n\in \mathbb{Z}$ since $R$ is right coherent, where $I^{+}=\underline{{\rm Hom}}_{\mathbb{Z}}(I,\overline{\mathbb{Q}/\mathbb{Z}})$. Thus the sequence $$\xymatrix@-1pc{\cdots\ar[r]^{} & {\rm Hom}_{C(R)}(F^{-1},I^{+}[-n]) \ar[r]& {\rm Hom}_{C(R)}(F^{0},I^{+}[-n])\ar[r] &{\rm Hom}_{C(R)}(F^{1},I^{+}[-n])\ar[r] &\cdots}$$ is exact for any $n\in \mathbb{Z}$, and so
$$\xymatrix@-1pc{\cdots\ar[r]^{} & \underline{{\rm Hom}}_{R}(F^{-1},I^{+}) \ar[r]& \underline{{\rm Hom}}_{R}(F^{0},I^{+})\ar[r] &\underline{{\rm Hom}}_{R}(F^{1},I^{+})\ar[r] &\cdots}$$ is exact. Hence the sequence $$\xymatrix@-1pc{\cdots\ar[r]^{} & I\overline{\otimes}_R F^1 \ar[r]& I\overline{\otimes}_R F^0\ar[r] &I\overline{\otimes}_R F^{-1}\ar[r] &\cdots}$$ is exact by \cite[Proposition 4.2.1(1)]{gr}. Therefore $C$ is Gorenstein flat.

$\Longleftarrow$) Suppose that $C$ is a Gorenstein flat complex. Then there exists an exact sequence $\xymatrix@-1pc{\cdots\ar[r]^{} & F^{1} \ar[r]& F^{0}\ar[r] &F^{-1}\ar[r] &\cdots}$ of flat complexes with $C\cong{\rm Im}(\xymatrix@-1pc{F^{0}\ar[r] &F^{-1})}$ and which remains exact after applying $I\overline{\otimes}_R-$ for any injective
$R^{\circ}$-complex $I$. Let $K\in \mathcal{\widetilde{F}\cap \text{dg}\widetilde{C}}$. Then we have a pure exact sequence $\xymatrix@-1pc{0\ar[r]^{} & K \ar[r]& K^{++}\ar[r] &K^{++}/K\ar[r] &0}$ by \cite[Proposition 5.1.4(4)]{gr}. Since $K\in \widetilde{\mathcal{F}}$, we get $K^{++}\in \widetilde{\mathcal{F}}$. So $K^{++}/K\in \widetilde{\mathcal{F}}$ by \cite[Lemma 4.7]{g04}. Thus the sequence $\xymatrix@-1pc{0\ar[r]^{} & K \ar[r]& K^{++}\ar[r] &K^{++}/K\ar[r] &0}$ is split. By \cite[Proposition 4.2.1(1)]{gr}, we have the following commutative diagram
$$\xymatrix@-1pc{ \cdots\ar[r] &( K^{+}\overline{\otimes}_R F^{-1})^{+}\ar[d]_{\cong} \ar[r]&(K^{+}\overline{\otimes}_R F^{0})^{+}  \ar[d]_{\cong}\ar[r] &
(K^{+}\overline{\otimes}_R F^{1})^{+}\ar[d]_{\cong} \ar[r] &\cdots \\
 \cdots\ar[r]&\underline{{\rm Hom}}_{R}(F^{-1},K^{++}) \ar[r]&\underline{{\rm Hom}}_{R}(F^{0},K^{++})  \ar[r] & \underline{{\rm Hom}}_{R}(F^{1},K^{++})\ar[r]&\cdots}$$ where the top row is exact since $K^{+}$ is injective. So the lower row is exact. Hence the sequence
 $$\xymatrix@-1pc{\cdots\ar[r]&\underline{{\rm Hom}}_{R}(F^{-1},K) \ar[r]&\underline{{\rm Hom}}_{R}(F^{0},K)  \ar[r] & \underline{{\rm Hom}}_{R}(F^{1},K)\ar[r]&\cdots}$$ is exact. In particular, the sequence$$\xymatrix@-1pc{\cdots\ar[r]^{} & {\rm Hom}_{C(R)}(F^{-1},K) \ar[r]& {\rm Hom}_{C(R)}(F^{0},K)\ar[r] &{\rm Hom}_{C(R)}(F^{1},K)\ar[r] &\cdots}$$ is exact. So $C$ is an \textbf{F}-Gorenstein flat complex.
\endproof

The following result will be used in the sequel.

\begin{Lem}\label{seq}  Let
$\xymatrix@-1pc{\cdots\ar[r]& X^{1} \ar[r]^{} & X^{0}\ar[r]^{} & X^{-1}\ar[r]^{} &\cdots}$ be a ${\rm Hom}_{C(R)}(-,\mathcal{\widetilde{A}\cap \text{dg}\widetilde{B}})$-exact sequence of complexes, then the sequence $\xymatrix@-1pc{\cdots\ar[r]^{} & X^{1}_n \ar[r]^{} & X^{0}_n\ar[r]^{} & X^{-1}_n\ar[r]^{} &\cdots}$ is
${\rm Hom}_R(-,\mathcal{A}\cap \mathcal{B})$-exact for any $n\in \mathbb{Z}$.\end{Lem}

\proof Let $K\in\mathcal{A}\cap \mathcal{B}$ and $n\in \mathbb{Z}$. Then $\overline{K}[n]\in \mathcal{\widetilde{A}\cap \text{dg}\widetilde{B}}$ by \cite[Lemma 3.4]{g04}. So we have the following exact sequence
$$\xymatrix@-1pc{\cdots\ar[r]&{\rm Hom}_{C(R)}(X^{-1}, \overline {K}[n]) \ar[r]_{ }     &{\rm Hom}_{C(R)}(X^{0}, \overline{K} [n])  \ar[r]_{ } & {\rm Hom}_{C(R)}(X_{1}, \overline{K}[n])  \ar[r]&\cdots.}$$
Using the standard adjunction of \cite[Lemma 3.1(2)]{g04}, we get the exact sequence
$$\xymatrix@-0.5pc{\cdots\ar[r]&\text{Hom}_R(X^{-1}_n,K)  \ar[r]_{ } &\text{Hom}_R(X^{0}_n, K)  \ar[r]_{ } & \text{Hom}_R(X^{1}_n, K)  \ar[r]&\cdots.}$$
This completes the proof.
\endproof

Now, we are in position to prove our main result, which gives a characterization of complexes in $\mathcal{G}(\widetilde{\mathcal{A}})$ and unifies \cite[Theorems 4.7]{hx} and \cite[Theorem 2.2]{yl}.

\begin{The}\label{T1} Let $C$ be an $R$-complex. Then $C\in \mathcal{G}(\widetilde{\mathcal{A}})$ if and only if $C_n\in\mathcal{G}(\mathcal{A})$ for any $n\in \mathbb{Z}$.
\end{The}

\proof
$\Longrightarrow)$ Assume that $C\in \mathcal{G}(\widetilde{\mathcal{A}})$. Then there exists a ${\rm Hom}_{C(R)}(-,\mathcal{\widetilde{A}\cap \text{dg}\widetilde{B}})$-exact exact sequence $\xymatrix@-1pc{\cdots\ar[r]^{} & A^{1} \ar[r]& A^{0}\ar[r] &A^{-1}\ar[r] &\cdots}$ with each $A^i\in \mathcal{\widetilde{A}}$ such that $C\cong{\rm Im}(\xymatrix@-1pc{A^{0}\ar[r] &A^{-1})}$. Now for any but fixed $n\in \mathbb{Z}$, by Lemma \ref{seq}, we have the following ${\rm Hom}_R(-,\mathcal{A}\cap \mathcal{B})$-exact exact sequence of modules in $\mathcal{A}$
$$\xymatrix{\cdots\ar[r]^{} & A^{1}_n \ar[r]^{} & A^{0}_n\ar[r]^{} & A^{-1}_n\ar[r]^{} & A^{-2}_n\ar[r]^{}&\cdots}$$ such that $C_n={\rm Im}(\xymatrix@-1pc{A^{0}_n\ar[r] &A^{-1}_n})$.
Hence $C_{n}\in\mathcal{G}(\mathcal{A})$.

$\Longleftarrow)$ Suppose that $C_{n}\in\mathcal{G}(\mathcal{A})$  for all $n\in \mathbb{Z}$. Then for any $n\in \mathbb{Z}$, there exists an exact sequence
$$\xymatrix{ 0\ar[r]^{} & C_{n} \ar[r]^{} & A_{n}\ar[r]^{} & L_{n}\ar[r]^{}&0,}$$ where $A_n\in \mathcal{A}$ and $L_{n}\in \mathcal{G}(\mathcal{A})$. These exact sequences induce a short exact sequence of complexes
$$\xymatrix{0\ar[r]^{} & \oplus_{n\in \mathbb{Z}}\overline{C_n}[n]\ar[r]^{} & \oplus_{n\in \mathbb{Z}}\overline{A_n}[n]\ar[r]^{} & \oplus_{n\in \mathbb{Z}}\overline{L_n}[n]\ar[r]^{} &0.}$$ Put $A^{-1}=\bigoplus_{n\in \mathbb{Z}}\overline{A_n}[n]$. It is easy to see that $A^{-1}\in \widetilde{\mathcal{A}}$. On the other hand, there is an obvious
(degreewise split) short exact sequence $$\xymatrix{0\ar[r]^{} & C\ar[r]^{\tiny\tiny\left(\begin{array}{c}
                                                                                1 \\
                                                                                \delta
                                                                              \end{array}\right)\;\;\;\;\;\;\;\;\;\;
} & \oplus_{n\in \mathbb{Z}}\overline{C_n}[n]\ar[r]^{\quad(-\delta,1)} & C[1] \ar[r]&0,}$$ where $\delta$ is the differential of $C$. Now let $\alpha:\xymatrix{C\ar[r]&A^{-1}}$ be the composite $$\xymatrix{C\ar[r]^{} & \oplus_{n\in \mathbb{Z}}\overline{C_n}[n]\ar[r]^{} & \oplus_{n\in \mathbb{Z}}\overline{A_n}[n].}$$ Then $\alpha$ is monoic since it is the composite of two monomorphisms. Denote ${\rm Coker}\alpha$ by $C^{-1}$. Then by Snake Lemma, we have a short exact sequence
$$\xymatrix{0\ar[r]^{} &C[1]\ar[r]^{} & C^{-1}\ar[r]^{} & \oplus_{n\in \mathbb{Z}}\overline{L_n}[n]\ar[r]^{} &0.}$$
Since each degree of $\oplus_{n\in \mathbb{Z}}\overline{L_n}[n]$ and $C[1]$ are in $\mathcal{G}(\mathcal{A})$, each degree of $C^{-1}$ belongs to $\mathcal{G}(\mathcal{A})$ by \cite[Proposition 3.3(1)]{yc}. Let $K\in \mathcal{\widetilde{A}\cap \text{dg}\widetilde{B}}$. Then $K\in \mathcal{\widetilde{A}\cap \widetilde{B}}=\widetilde{\mathcal{A}\cap \mathcal{B}}$ by \cite[Theorem 3.12]{g04}. Thus $K\cong\prod_{n\in \mathbb{Z}}\overline{{\rm Z}_n(K)}[n]$ by \cite[Lemma 4.1]{ldy}. Hence $${\rm Ext}_{C(R)}^1(C^{-1},K)\cong\prod_{n\in \mathbb{Z}}{\rm Ext}_{C(R)}^1(C^{-1},\overline{{\rm Z}_n(K)}[n])\cong\prod_{n\in \mathbb{Z}}{\rm Ext}_{R}^1(C^{-1}_n,{\rm Z}_n(K))=0,$$ where the second isomorphism follows from \cite[Lemma 3.1(2)]{g04} and the last equality follows from Remark \ref{rem1}(1). This implies that $\xymatrix@-1pc{0\ar[r]^{} & C\ar[r]^{} & A^{-1}\ar[r]^{} & C^{-1}\ar[r]^{} &0}$
is ${\rm Hom}_{C(R)}(-,\mathcal{\widetilde{A}\cap \text{dg}\widetilde{B}})$-exact.
Notice that $C^{-1}$ has the same property as $C$, so we can use the same procedure to construct a ${\rm Hom}_{C(R)}(-,\mathcal{\widetilde{A}\cap \text{dg}\widetilde{B}})$-exact exact sequence of complexes
$$\xymatrix{0\ar[r]& C\ar[r]^{}&A^{-1}\ar[r]^{}&A^{-2}\ar[r]& \cdots,} \eqno(\dagger)\vspace*{-2mm}$$
where each $A^i$ is an $\mathcal{A}$-complex.

Since $(\mathcal{\widetilde{A}, \text{dg}\widetilde{B}})$ is a complete cotorsion pair, we have a short exact sequence $\xymatrix@-1pc{0\ar[r]^{} & C^1\ar[r]^{} & A^{0}\ar[r]^{} & C\ar[r]^{} &0}$, where $A^0\in \widetilde{\mathcal{A}}$ and $C^1\in\text{dg}\widetilde{\mathcal{B}}$. Note that $C_n\in \mathcal{G}(\mathcal{A})$ for any $n\in Z$, this sequence is ${\rm Hom}_{C(R)}(-,\mathcal{\widetilde{A}\cap \text{dg}\widetilde{B}})$-exact by a similarly discussion as above. Also, it follows from the exact sequence and \cite[Proposition 3.3(1)]{yc} that each $C_n^1\in \mathcal{G}(\mathcal{A})$ for any $n\in \mathbb{Z}$. Thus we can continuously use the same method to construct a ${\rm Hom}_{C(R)}(-,\mathcal{\widetilde{A}\cap \text{dg}\widetilde{B}})$-exact exact sequence
 $$\xymatrix{\cdots\ar[r]& A^1\ar[r]^{}&A^0\ar[r]^{}&C\ar[r]& 0,} \eqno(\ddagger)\vspace*{-2mm}$$
where each $A^i$ is an $\mathcal{A}$-complex.

Finally, gluing the sequences $(\dagger)$ and $(\ddagger)$ together, one has a ${\rm Hom}_{C(R)}(-,\mathcal{\widetilde{A}\cap \text{dg}\widetilde{B}})$-exact exact sequence of complexes
$$\xymatrix{\cdots\ar[r]& A^{1} \ar[r]^{} & A^{0}\ar[r]^{} &A^{-1}\ar[r]^{}& A^{- 2}\ar[r]^{} &\cdots}$$ with all $A^{i}\in \widetilde{\mathcal{A}}$ such that $C\cong{\rm Im}(\xymatrix@-1pc{A^{0}\ar[r] &A^{-1}})$. Hence $C\in \mathcal{G}(\widetilde{\mathcal{A}})$.
\endproof

Let $\mathcal{D}$ be an abelian category with enough projective objects and injective objectives. Recall that a class $\mathcal{X}$ of objects of $\mathcal{D}$ is said to be projectively resolving (resp., injectively resolving) if it is closed under extensions and
kernels of surjections (resp., cokernels of injections), and it contains all projective (resp., injective) objects of $\mathcal{D}$.

\begin{Cor}\label{}  $\mathcal{G}(\widetilde{\mathcal{A}})$ is projectively resolving.\end{Cor}

\proof Clearly, $\widetilde{\mathcal{P}}\subseteq\mathcal{\widetilde{A}}\subseteq\mathcal{G}(\widetilde{\mathcal{A}})$. Let $\xymatrix@-1pc{0\ar[r]&C'\ar[r]&C\ar[r]&C''\ar[r]&0}$ be a short exact sequence of complexes with $C''\in\mathcal{G}(\widetilde{\mathcal{A}})$. Then for any $n\in \mathbb{Z}$, in the exact sequence $\xymatrix@-1pc{0\ar[r]&C'_n\ar[r]&C_n\ar[r]&C''_n\ar[r]&0,}$ $C''_n\in \mathcal{G}(\mathcal{A})$ by Theorem \ref{T1}. So $C'_n\in \mathcal{G}(\mathcal{A})$ if and only if $C_n\in \mathcal{G}(\mathcal{A})$ by \cite[Proposition 3.3(1)]{yc}. Hence $C'\in\mathcal{G}(\widetilde{\mathcal{A}})$ if and only if $C\in\mathcal{G}(\widetilde{\mathcal{A}})$ by Theorem \ref{T1}.
Now the result follows.
\endproof

\begin{Cor}\label{third} Let $\xymatrix@-1pc{0\ar[r] &C'\ar[r] &C\ar[r] &C''\ar[r]&0}$ be an exact sequence of complexes.
If $C',C$ are belong to $\mathcal{G}(\widetilde{\mathcal{A}})$, then $C''\in\mathcal{G}(\widetilde{\mathcal{A}})$
if and only if ${\rm Ext}^{1}_{C(R)}(C'',K)=0$ for any $K\in \widetilde{\mathcal{A}}\cap\text{dg}\widetilde{\mathcal{B}}$.\end{Cor}

\proof $\Longrightarrow$) It is obvious.

$\Longleftarrow$)  Let $n\in \mathbb{Z}$. Consider the exact sequence of $R$-modules
$$\xymatrix{0\ar[r]& C'_{n}\ar[r]&C_{n}\ar[r]&C''_{n}\ar[r]&0.}$$
By Theorem \ref{T1}, $C'_{n},C_{n}$ are belong to $\mathcal{G}(\mathcal{A})$. Let $K\in \mathcal{A}\cap\mathcal{B}$. Then $\overline{K}[n]\in \widetilde{\mathcal{A}}\cap\text{dg}\widetilde{\mathcal{B}}$. Thus ${\rm Ext}^{1}_R(C''_{n},K)\cong{\rm Ext}^{1}_{\mathcal{C}(R)}(C'',\overline{K}[n])$=0 by  \cite[Lemma 3.1(2)]{g04} and the hypothesis.
Hence $C''_{n}\in \mathcal{G}(\mathcal{A})$ by \cite[Proposition 3.3(2)]{yc}. Therefore $C''\in\mathcal{G}(\widetilde{\mathcal{A}})$ by Theorem \ref{T1}.\endproof

By Proposition \ref{prop}, Theorem \ref{T1} and \cite[Lemma 3.2]{hx}, we immediately get that

\begin{Cor}(\cite[Theorem 3.1]{yl}) Let $C$ be an $R$-complex. If $R$ is a right coherent ring, then $C$ is Gorenstein flat if and only if $C_n$ is a Gorenstein flat $R$-module for any $n\in \mathbb{Z}$.\end{Cor}

\section{Stability of Gorenstein categories with respect to cotorsion pairs}\label{sec1}
The stability of Gorenstein categories was initiated by Sather-Wagstaff, Sharif and White \cite{ssw}. They proved that if $R$ is a commutative ring, then an $R$-module $M$ is a Gorenstein projective (resp., injective) module if and only if there exists an exact sequence of Gorenstein projective (resp., injective) $R$-modules
$\xymatrix@-1pc{G=\cdots \ar[r]&G_{1} \ar[r]^{\delta_{1}}  &G_{0}  \ar[r]^{\delta_{0}} & G_{-1}  \ar[r]&\cdots}$
such that the complexes ${\rm Hom}_R(H,G)$ and ${\rm Hom}_R(G,H)$
are exact for each Gorenstein projective (resp., injective) $R$-module $H$ and $M={\rm Im}\delta_{0}$. This was developed by Bouchiba \cite{b}, Xu and Ding \cite{xd}, respectively. They showed, via different methods, that over any ring $R$, an $R$-module $M$ is Gorenstein projective (resp., injective) if and only if there exists an exact sequence of Gorenstein projective (resp., injective) $R$-modules
$\xymatrix@-1pc{G=\cdots \ar[r]&G_{1} \ar[r]^{\delta_{1}}  &G_{0}  \ar[r]^{\delta_{0}} & G_{-1}  \ar[r]&\cdots}$
such that the complex ${\rm Hom}_R(G,H)$ (resp., ${\rm Hom}_R(H,G)$)
is exact for any projective (resp., injective) $R$-module $H$ and $M={\rm Im}\delta_{0}$. For more details, see \cite{b}. The stabiltity of Gorenstein flat  $R$-module  has been treated by Bouchiba and Khaloui \cite{bk}, Xu and Ding \cite{xd}, Yang and Liu \cite{yl12}, respectively. By using totally different techniques, they showed that over a left GF-closed ring $R$ (a ring $R$ over which the class of the Gorenstein
flat $R$-modules is closed under extensions), an $R$-module $M$ is Gorenstein flat if and only if there exists an exact sequence of Gorenstein flat $R$-modules $\xymatrix@-1pc{G=\cdots \ar[r]&G^{1} \ar[r]^{\delta^{1}}  &G^{0}  \ar[r]^{\delta^{0}} & G^{-1}  \ar[r]&\cdots}$ such that the complex $I\otimes_R G$ is exact
for each Gorenstein injective (or injective) $R^{\circ}$-module $I$ and $M={\rm Im}\delta^{0}$. By using Theorem \ref{T1}, in this section, we investigate the stability of  $\mathcal{G}(\mathcal{A})$ and $\mathcal{G}(\mathcal{\widetilde{A}})$.

The next result shows that the category $\mathcal{G}(\mathcal{A})$ possesses of stability, which is a generalization of \cite[Theorem A]{ssw}, \cite[Theorem A]{xd} and \cite[Theorem 3.8]{hx}.

\begin{The}\label{T2} Let $M$ be an $R$-module. Then the following statements are equivalent:
\par(1) $M\in\mathcal{G}(\mathcal{A}).$
\par(2) there exists a both ${\rm Hom}_{R}(\mathcal{G}(\mathcal{A}),-)$-exact and ${\rm Hom}_{R}(-,\mathcal{G}(\mathcal{A}))$-exact exact sequence $\xymatrix@-1pc{\cdots\ar[r]^{} & G_{1} \ar[r]& G_{0}\ar[r] &G_{-1}\ar[r] &\cdots}$ in $\mathcal{G}(\mathcal{A})$ such that $M\cong{\rm Im}(\xymatrix@-1pc{G_{0}\ar[r] &G_{-1}).}$
\par(3) there exists a ${\rm Hom}_{R}(-,\mathcal{G}(\mathcal{A}))$-exact exact sequence $\xymatrix@-1pc{\cdots\ar[r]^{} & G_{1} \ar[r]& G_{0}\ar[r] &G_{-1}\ar[r] &\cdots}$
in $\mathcal{G}(\mathcal{A})$ such that $M\cong{\rm Im}(\xymatrix@-1pc{G_{0}\ar[r] &G_{-1}).}$
\par(4) there exists a ${\rm Hom}_{R}(-,\mathcal{A})$-exact exact sequence $\xymatrix@-1pc{\cdots\ar[r]^{} & G_{1} \ar[r]& G_{0}\ar[r] &G_{-1}\ar[r] &\cdots}$
in $\mathcal{G}(\mathcal{A})$ such that $M\cong{\rm Im}(\xymatrix@-1pc{G_{0}\ar[r] &G_{-1}).}$
\par(5) there exists a ${\rm Hom}_{R}(-,\mathcal{A}\cap \mathcal{B})$-exact exact sequence $\xymatrix@-1pc{\cdots\ar[r]^{} & G_{1} \ar[r]& G_{0}\ar[r] &G_{-1}\ar[r] &\cdots}$
in $\mathcal{G}(\mathcal{A})$ such that $M\cong{\rm Im}(\xymatrix@-1pc{G_{0}\ar[r] &G_{-1}).}$
\end{The}

\proof (1)$\Longrightarrow$(2)$\Longrightarrow$(3)$\Longrightarrow$(4)$\Longrightarrow$(5) are clear.

(5)$\Longrightarrow$(1) Assume that there is a ${\rm Hom}_{R}(-,\mathcal{A}\cap \mathcal{B})$-exact exact sequence
$$\xymatrix{G=\cdots\ar[r]^{} & G_{1} \ar[r]& G_{0}\ar[r] &G_{-1}\ar[r] &G_{-2}\ar[r]^{}&\cdots}$$
in $\mathcal{G}(\mathcal{A})$ such that $M\cong{\rm Z}_{-1}(G)$. Then $G\in\mathcal{G}(\widetilde{\mathcal{A}})$ by Theorem \ref{T1}. Thus there exists a ${\rm Hom}_{C(R)}(-,\mathcal{\widetilde{A}\cap \text{dg}\widetilde{B}})$-exact exact sequence $$\xymatrix{\cdots\ar[r]^{} & A^{1} \ar[r]^{\sigma^{1}}& A^{0}\ar[r]^{\sigma^{0}} &A^{-1}\ar[r]^{\sigma^{-1}}&A^{-2}\ar[r] &\cdots}$$ with each $A^i\in \mathcal{\widetilde{A}}$ such that $G\cong{\rm Ker}\sigma^{-1}$. Set $K^i={\rm Ker}\sigma^{i}$ for $i\in \mathbb{Z}$. Then $K^i\in\mathcal{G}(\widetilde{\mathcal{A}})$ and $K^i$ is exact for any $i\in \mathbb{Z}$ since $K^{-1}=G$ and all $A^i$ are exact. So, by \cite[Lemma 4.15(1)]{ldy}, we have the following exact sequence \vspace*{-2mm}$$\xymatrix{\cdots\ar[r]^{} & {\rm Z}_{-1}(A^{1}) \ar[r]^{{\rm Z}_{-1}(\sigma^1)} &{\rm Z}_{-1}(A^{0})\ar[r]^{{\rm Z}_{-1}(\sigma^0)} &{\rm Z}_{-1}(A^{-1})\ar[r]^{} &\cdots}\eqno(\natural)\vspace*{-2mm}$$ with each ${\rm Z}_{-1}(A^i)\in \mathcal{A}$, such that $M\cong{\rm Z}_{-1}(G)={\rm Ker}({\rm Z}_{-1}(\sigma^{-1}))$. To show $M\in\mathcal{G}(\mathcal{A})$, we need only to show that the sequence $(\natural)$ is ${\rm Hom}_R(-,\mathcal{A}\cap \mathcal{B})$-exact.

Let $H\in\mathcal{A}\cap \mathcal{B}$, it suffices to show that ${\rm Ext}_R^1({\rm Z}_{-1}(K^i),H)=0$ for all $i\in \mathbb{Z}$. Since each $K^i\in\mathcal{G}(\widetilde{\mathcal{A}})$, all $K^i_n\in\mathcal{G}(\mathcal{A})$ by Theorem \ref{T1}.
Thus, for any $i\in \mathbb{Z}$, the sequence
$$\xymatrix{ 0\ar[r]_{ }       &{\rm Hom}_R(K^{i-1},H)\ar[r]_{ }&{\rm Hom}_R(A^{i},H) \ar[r]_{ } &{\rm Hom}_R(K^i,H)\ar[r]_{ }&0 }$$ is exact. By the hypothesis, ${\rm Hom}_R(K^{-1},H)$ is exact. Note that ${\rm Hom}_R(A^i,H)$ is exact for each $i\in \mathbb{Z}$, then ${\rm Hom}_R(K^i,H)$ is exact for any $i\in \mathbb{Z}$. Hence
${\rm Ext}^1_R({\rm Z}_{-1}(K^i),H)=0$ since each $K_{0}^i\in \mathcal{G}(\mathcal{A})$. Thus the sequence $(\natural)$ is ${\rm Hom}_R(-,\mathcal{A}\cap \mathcal{B})$-exact, as desired.
\endproof

Finally, by applying Theorem \ref{T1} and Theorem \ref{T2}, we can achieve the following stability result for $\mathcal{G}(\widetilde{\mathcal{A}})$, which is a unification of \cite[Theorem 3.1]{xd} and \cite[Theorem 4.11]{hx}.

\begin{The}\label{T3} Let $C$ be a complex of $R$-modules. Then the following statements are equivalent:
\par(1) $C\in\mathcal{G}(\widetilde{\mathcal{A}}).$
\par(2) there exists a both ${\rm Hom}_{C(R)}(\mathcal{G}(\widetilde{\mathcal{A}}),-)$-exact and ${\rm Hom}_{C(R)}(-,\mathcal{G}(\widetilde{\mathcal{A}}))$-exact exact sequence $\xymatrix@-1pc{\cdots\ar[r]^{} & G^{1} \ar[r]& G^{0}\ar[r] &G^{-1}\ar[r] &\cdots}$
in $\mathcal{G}(\widetilde{\mathcal{A}})$ such that $C\cong{\rm Im}(\xymatrix@-1pc{G^{0}\ar[r] &G^{-1}).}$
\par(3) there is a ${\rm Hom}_{C(R)}(-,\mathcal{G}(\widetilde{\mathcal{A}}))$-exact exact sequence $\xymatrix@-1pc{\cdots\ar[r]^{} & G^{1} \ar[r]& G^{0}\ar[r] &G^{-1}\ar[r] &\cdots}$
in $\mathcal{G}(\widetilde{\mathcal{A}})$ such that $C\cong{\rm Im}(\xymatrix@-1pc{G^{0}\ar[r] &G^{-1}).}$
\par(4) there is a ${\rm Hom}_{C(R)}(-,\widetilde{\mathcal{A}})$-exact exact sequence $\xymatrix@-1pc{\cdots\ar[r]^{} & G^{1} \ar[r]& G^{0}\ar[r] &G^{-1}\ar[r] &\cdots}$
in $\mathcal{G}(\widetilde{\mathcal{A}})$ such that $C\cong{\rm Im}(\xymatrix@-1pc{G^{0}\ar[r] &G^{-1}).}$
\par(5) there is a ${\rm Hom}_{C(R)}(-,\widetilde{\mathcal{A}}\cap {\rm dg}\widetilde{\mathcal{B}})$-exact exact sequence $\xymatrix@-1pc{\cdots\ar[r]^{} & G^{1} \ar[r]& G^{0}\ar[r] &G^{-1}\ar[r] &\cdots}$
in $\mathcal{G}(\widetilde{\mathcal{A}})$ such that $C\cong{\rm Im}(\xymatrix@-1pc{G^{0}\ar[r] &G^{-1}).}$
\end{The}

\proof (1)$\Longrightarrow$(2)$\Longrightarrow$(3)$\Longrightarrow$(4)$\Longrightarrow$(5) are trivial.

(5)$\Longrightarrow$(1) Suppose that there exists a ${\rm Hom}_{C(R)}(-,\widetilde{\mathcal{A}}\cap {\rm dg}\widetilde{\mathcal{B}})$-exact exact sequence $$\xymatrix{\cdots\ar[r] & G^{1} \ar[r]^{\sigma^1}& G^{0}\ar[r]^{\sigma^0} &G^{-1}\ar[r]^{\sigma^{-1}} &G^{-2}\ar[r]^{}&\cdots}$$
in $\mathcal{G}(\widetilde{\mathcal{A}})$ such that $C\cong{\rm Im}{\sigma^0}$. Then for any $n\in \mathbb{Z}$, by Lemma \ref{seq}, we have the following ${\rm Hom}_R(-,\mathcal{A}\cap \mathcal{B})$-exact exact sequence of modules
$$\xymatrix{\cdots\ar[r]^{} & G^{1}_n \ar[r]^{\sigma^{1}_n} & G^{0}_n\ar[r]^{\sigma^{0}_n} & G^{-1}_n\ar[r]^{\sigma^{-1}_n} & G^{-2}_n\ar[r]^{}&\cdots}$$ such that $C_n\cong{\rm Im}\sigma_0^n$. By Theorem \ref{T1}, $G^{i}_n\in \mathcal{G}(\mathcal{A})$ for each $i\in \mathbb{Z}$. Thus $C_n\in \mathcal{G}(\mathcal{A})$ by Theorem \ref{T2}.
Hence $C\in\mathcal{G}(\widetilde{\mathcal{A}})$ by Theorem \ref{T1}.
\endproof

\vspace{0.5cm}
{\small
{\em Authors' addresses}:
{\em Renyu Zhao}(corresponding author), Department of Mathematics, Northwest Normal University, Anning East road No. 967,
Lanzhou, Gansu, 730070, P.R.China, and Department of Mathematics, Nanjing University, Nanjing 210093, P.R. China, e-mail: \texttt{zhaory@\allowbreak nwnu.edu.cn};
{\em Pengju Ma}, Department of Mathematics, Northwest Normal University, Anning East road No. 967,
Lanzhou, Gansu, 730070, P.R.China,
 e-mail: \texttt{2642293920@\allowbreak qq.com}.

}

\end{document}